\documentclass[reqno,10pt]{amsart}
\usepackage{amssymb,amsmath,amsthm,amsfonts,color}

\usepackage{mathrsfs,dsfont, comment,mathscinet}

\usepackage{enumerate,esint}
\usepackage[left=2.8cm,right=2.8cm,top=2.8cm,bottom=2.8cm]{geometry}
\usepackage{mathtools}
\usepackage{graphicx}
\usepackage[colorinlistoftodos]{todonotes}
\usepackage{epstopdf}
\usepackage[colorlinks=true, pdfstartview=FitV, linkcolor=blue, 
            citecolor=blue, urlcolor=blue]{hyperref}

\allowdisplaybreaks
\usepackage{verbatim}
\usepackage{pgfplots}
\pgfplotsset{width=7cm,compat=1.8}
\pgfplotsset{
    colormap/outside/.style={
        colormap=
            {outside}{
            rgb255(0cm)=(110,110,255);
            rgb255(1cm)=(20,20,255);
            }
    },
    colormap/outside,
    colormap/inside/.style={
        colormap={inside}{
            rgb255(0cm)=(20,20,255);
            rgb255(1cm)=(220,220,255);
        }
    },
    colormap/inside
}

\numberwithin{equation}{section}
\theoremstyle{plain}
\newtheorem{theorem}{Theorem}[section]

\newtheorem{lemma}[theorem]{Lemma}

\theoremstyle{definition}
\newtheorem{definition}[theorem]{Definition}

\newcommand\R{\mathbb R}

\newcommand\mthree{\mathbb{M}^{3\times 3}}
\newcommand\nh{\nabla_h}

\newcommand\wk{\rightharpoonup}

\newcommand\be[1]{\begin{equation}\label{#1}}
\newcommand\ee{\end{equation}}
\newcommand\ba[1]{\begin{align}\label{#1}}
\newcommand\bas{\begin{align*}}
\newcommand\nn{\nonumber}
\newcommand\ep{\varepsilon}
\renewcommand{\O}{\Omega}
\newcommand{\md}{{\rm d}}

\newcommand\UUU{\color{black}}

\newcommand\EEE{\color{black}}
\newcommand\RRR{\color{black}}

\newcommand{\MMM}{\color{black}}

\title [{Magnetoelastic thin films}]{Magnetoelastic thin films \UUU at
  large strains}%\EEE under injectivity constraints}
\author{Elisa Davoli}
\address[Elisa Davoli]{Institute of Analysis and Scientific Computing, Vienna University of Technology,
   Wiedner Hauptstrasse 8-10, 1040 Vienna, Austria} 
\email{elisa.davoli@univie.ac.at}

\author{Martin Kru\v{z}\'{i}k}
\address[Martin Kru\v{z}\'{i}k]{ Czech Academy of Sciences, Institute of Information Theory
   and Automation, Pod Vod\'{a}renskou v\v{e}\v{z}\'{i}
    4, 180 00 Prague, Czechia}
 \email{kruzik@utia.cas.cz}

\author{Paolo Piovano}
\address[Paolo Piovano]{Faculty of Mathematics, University of Vienna,
   Oskar-Morgenstern-Platz 1, A-1090 Vienna, Austria}
\email{paolo.piovano@univie.ac.at}

\author{Ulisse Stefanelli} 
\address[Ulisse Stefanelli]{Faculty of Mathematics, University of
  Vienna, Oskar-Morgenstern-Platz 1, A-1090 Vienna, Austria,
Vienna Research Platform on Accelerating
  Photoreaction Discovery, University of Vienna, W\"ahringerstra\ss e 17, 1090 Wien, Austria,
 \& Istituto di
  Matematica Applicata e Tecnologie Informatiche {\it E. Magenes}, via
  Ferrata 1, I-27100 Pavia, Italy
}
\email{ulisse.stefanelli@univie.ac.at}

\EEE

\begin{document} 
\vskip .2truecm
\begin{abstract}
Starting from \UUU the three-dimensional setting, \EEE  
we derive a limit model of a thin \UUU magnetoelastic \EEE film by
means of $\Gamma$-convergence \UUU techniques. \EEE As magnetization vectors are \UUU
defined on the \EEE  elastically deformed configuration, our model
\UUU features both Lagrangian \EEE and Eulerian \UUU terms. \EEE This
calls for qualifying admissible three-dimensional deformations of
planar domains \UUU in terms of injectivity.  
In addition, a careful treatment of the Maxwell system in the deformed
film is required.  \EEE
\end{abstract}
\maketitle

\section{Introduction}
 Magnetoelasticity describes the mechanical behavior of solids under
magnetic effects. % and  viceversa.   
The magnetoelastic
coupling is based on the presence of small magnetic domains in the
material \cite{hubert-schaefer}.
In the absence of an external magnetic field,  these magnetic domains are randomly
 oriented. \UUU When \EEE an external magnetic field \UUU is applied, \EEE 
 \UUU the \UUU mesostructure  of magnetic domains changes by magnetic-domain wall
motion,  \RRR by \UUU magnetization-vector rotation, and, for some specific alloys, by magnetic-field-driven
martensitic-variant transformation. The net effect is a magnetically
induced deformation in the body. Conversely,
mechanical deformations modify the magnetic response of a specimen by
influencing the magnetic anisotropy of the domains, so that the
magnetic and the mechanical behavior of the material are fully
coupled. \EEE
  We refer \RRR to, e.g., \EEE \cite{brown,desim,desimone-james,jam-kin}  for  an exposition  on the 
\UUU foundations \EEE of magnetoelasticity \UUU and to
\cite{kruzik-prohl} for some related mathematical considerations. \EEE

 The mathematical modeling  of magnetoelasticity is a lively area
of research,  triggered by the interest in the  so-called {\it multifunctional} 
materials. Among these one has to mention rare-earth alloys
such as TerFeNOL and GalFeNOL,  as well as ferromagnetic shape-memory
alloys as Ni$_2$MnGa, NiMnInCo,
NiFeGaCo, FePt, FePd,
among others \cite{James-Wuttig}.  These materials exhibit \UUU a
remarkable \EEE
magnetostrictive behavior, \UUU for \EEE reversible strains as large as 10\% can be activated
by the imposition of relatively moderate magnetic fields. This strong
magnetoelastic coupling makes them relevant in a wealth of innovative applications including
sensors and actuators \UUU \cite{davino}. \EEE  

\UUU The aim of this paper is to present a   model of a thin film \EEE
 undergoing large strain deformations in the membrane
regime. \UUU This will be inferred from a variational
dimension-reduction procedure from a corresponding three-dimensional
model at large strains. \EEE

Dimension-reduction \UUU techniques \EEE play an important  role in
nonlinear analysis and numeric\MMM s, \EEE \UUU for \EEE they allow simpler
computational approaches, still  preserving  the  main features of the
\UUU corresponding \EEE bulk model. 
% Therefore, it is of  great importance  that  a clear relationship
% between the full three-dimensional problem and its lower-dimensional
% counterpart  is  made rigorous.  
The last decades have witnessed remarkable \UUU progresses on
dimension reduction by \EEE variational methods, particularly by
$\Gamma$-convergence \cite{braides,DalMaso:93}, together with
quantitative  rigidity estimates \cite{FrieseckeJamesMueller:02}.
Among the \UUU many results on the elastic response of low-dimensional
objects, \EEE we mention the rigorous justification of membrane theory \cite{Dret1, Dret2}, bending theory \cite{FrieseckeJamesMueller:02, Pantz}, and von K\'{a}rm\'{a}n theory \cite{hierarchy, lecumberry} for plates as variational limits of nonlinear three-dimensional elasticity for vanishing thickness. In particular, we refer to \cite{hierarchy} for the derivation of a hierarchy of different plate models and for a thorough literature review.  
 
 A  rigorous  derivation of \UUU a model for magnetic \EEE thin films
 has been first obtained  in \cite{gioia-james}. A  rate-independent
 evolution of Kirchhoff-Love magnetic plates together with the passage
 from three-dimensional linearized magnetoelasticity  to the
 corresponding two-dimensional theory is the subject of
 \cite{kruzik.stefanelli.zanini}. \UUU Magnetostriction in thin \RRR films \EEE has
 been considered, also from the numerical viewpoint, in \cite{Liakhova,Luskin06,Luskin07}. 
 With respect to \UUU these results, this paper presents a fundamental novelty \RRR as it represents the first rigorous analytical treatment including also   \UUU the {\it large-strain} magnetoelastic regime. \EEE  
 %With respect to \UUU these results, this paper presents a fundamental novelty, \EEE for  we  \UUU focus on the {\it large-strain} magnetoelastic regime. \EEE  

  In \RRR the \EEE classical dimension reduction for {\it small-strain} elastic thin plates, the
  analysis is set in cylindrical domains whose height\MMM s \EEE depend on a
  thickness parameter eventually tending to zero. The same \UUU
  setting \EEE applies in magnetoelasticity. \UUU Under the small-deformations assumption, the magnetization may be assumed to be directly
  defined on the reference configuration. This simplification is
  however not amenable in the large-strain regime, for the
  magnetization is defined on the {\it deformed} configuration
  instead. The latter is however a priori not known, as it depends on
  the deformation itself. \EEE  In particular, this naturally leads to
  a mixed \UUU Eulerian-Lagrangian \EEE formulation of the problem.  \UUU Compared with previous small-strain
  contributions, the mathematical framework of this work  is hence \EEE much
  more involved. \UUU A distinctive difficulty arises from the need of
  ensuring \EEE  that admissible deformations are globally
  injective. In the bulk, this can be \UUU achieved \EEE  by imposing
  the so-called Ciarlet-Ne\v{c}as condition
  \cite{ciarlet-necas}. \UUU For
  films, however, no comparable \EEE condition, i.e., allowing for a
  variational approach, \UUU seems to be available. \EEE  A further
  difficulty is represented by \UUU the Maxwell system, which is
  formulated in actual space. \EEE In order to \UUU identify \EEE the
  asymptotic behavior of the stray field,  we have to characterize the
  limiting differential constraints in weak form \UUU by \EEE keeping track of the deformed configuration. 

 \UUU The main result of the paper \EEE is the derivation of \UUU a
 variational model for \RRR  thin-film specimens \EEE as a $\Gamma$-limit
 of  a suitably scaled energies of a bulk model for vanishing
 thickness.  In Theorem~\ref{prop:first-liminf} we prove in full
 generality the $\Gamma$-$\liminf$ inequality, showing that our limit
 energy functional always represents a lower bound for the asymptotic
 behavior of the three-dimensional energy functionals. If the limit
 film deformation is {\it approximately injective} in the sense of
 Definition~\ref{AID}, we show that  the  $\Gamma$-$\liminf$ is indeed
 the largest lower semicontinuous lower bound  for the
 magnetoelastic-plate functionals as the thickness goes to zero, \EEE
 i.e., it is the $\Gamma$-limit; cf.~ Theorem~\ref{prop:limsup}.
 Here, the approximate injectivity means that there is a sequence of
 deformations of the bulk which are globally injective and converge in
 a suitable sense to the film deformation.  Additionally, in
 Theorem~\ref{thm:complete-gamma-conv}  we prove a complete \EEE
 $\Gamma$-convergence result  under the additional assumption that the
 admissible three-dimensional deformations satisfy a suitable
 injectivity requirement \UUU which guarantees \EEE that  the limit deformation of the film is globally injective.

 The paper is organized as follows. In Section \ref{sec:setting} we
 introduce the mathematical setting of the problem. Section
 \ref{sec:gamma} \UUU is devoted \MMM to the statements \UUU of all results,  
 Section \ref{sec:proofs} contains all  proofs. %\UUU the proof of \EEE Theorem~\ref{thm:complete-gamma-conv}.

\section{Setting of the problem}
\label{sec:setting}

%\subsection{Notation and preliminaries}\label{subsec:notation}
We use the standard notation for
Sobolev and Lebesgue spaces, i.e., \UUU $W^{k,p} $ and
$L^p $ \EEE \cite{adams}.  If  $A\in\R^{3\times 2}$ and $b\in\R^3$
we write $(A|b)\in\R^{3\times 3}$  for a matrix whose first two
columns are created by the first two columns of $A$ and the third one
by the vector $b$. \EEE  The set of proper rotations is denoted by
SO$(3):=\{R\in\R^{3\times 3}:\, R^\top R= RR^\top=\RRR \text{Id}\EEE, \ {\rm det}\, R=1\}$ where \RRR \text{Id}\EEE is the identity matrix.

\UUU Let \EEE $\omega \subset\R^2$ \UUU be \EEE  a bounded Lipschitz
domain representing the planar reference configuration  of the film,
\UUU define the reference configuration of a thin magnetoelastic plate as 
$$\Omega_h:=\omega\times\Big(-\frac{h}{2},\frac{h}{2}\Big),$$
  and  set
$\O:=\O_1$. \EEE In the
expression above, $h>0$ represents the thickness of the plate,
eventually bound to go to zero. \UUU Correspondingly, we will consider
limits \RRR as $h \to 0$ of \UUU sequences of functionals by  means of 
$\Gamma$-convergence \cite{DalMaso:93}. This is a standard approach to \MMM characterize 
\UUU the
limiting \EEE behavior of a sequence of bulk \UUU energies for
specimens of very small thickness. \EEE 
 %in one direction vanishes. \EEE        
Assume that $X$ is a subset of a reflexive Banach space. We say that $\{I_h\}_{h>0}$ for  $I_h:X\to\R\cup\{\infty\}$ $\Gamma$-converges to $I:X\to\R\cup\{+\infty\}$   if the following conditions hold simultaneously:
\begin{subequations}\label{Gamma:def}
\begin{align}
& \UUU \zeta_h\stackrel{X}{\to} \zeta  \ \Rightarrow \ \liminf_{h\to 0}I_h(\zeta_h)\ge I(\zeta),
\\
&\UUU \forall\, \zeta \in X\ \exists \{\hat\zeta_h\}_{h>0}\subset X: \quad
\hat\zeta_h\stackrel{X}{\to} \zeta \ \text{and} \ 
   \limsup_{h\to 0}I_h(\hat\zeta_h)=I(\zeta),
\end{align}
\end{subequations}
\UUU where the symbol $\stackrel{X}{\to} $ indicates \RRR the convergence with respect to a properly chosen (weak) topology in $X$. \EEE 
If \eqref{Gamma:def} holds we say that $I$ is the $\Gamma$-limit of
$\{I_h\}_{h>0}$ \UUU (with respect to that topology). 

 \EEE

% Let $\omega\subset\R^2$ be a bounded Lipschitz domain, and assume that the set
% $$\Omega_h:=\omega\times\Big(-\frac{h}{2},\frac{h}{2}\Big)$$
% is the reference configuration of a thin magnetoelastic plate. In
% the expression above, $h>0$ represents the thickness of the plate,
% eventually bound to go to zero.
\UUU The state of the magnetoelastic material is defined in terms of
its {\it deformation} $w$ and its magnetization
$m$. The deformation $w: \Omega_h \to \R^3$ is required to belong to  $W^{1,p}(\Omega_h;\R^3)$
for some given $$p>3,$$ to be orientation-preserving, namely, $\det \nabla w >0$
almost everywhere, and to \EEE satisfy   the Ciarlet-Ne\v{c}as
condition \cite{ciarlet-necas}  
\be{eq:ciarlet-necas}
\int_{\Omega_h} {\rm det }\,\nabla w \,dx\leq
\mathcal{L}^3(w(\Omega_h))
\ee
\UUU where $\mathcal{L}^3$ stands for the three-dimensional Lebesgue
measure.  In particular, $w$ \RRR is \UUU identified with the unique continuous
representative in the equivalence class. The magnetization $m$ \MMM is set \EEE on the open deformed configuration,
namely, $m:\Omega^w_h \to {\mathbb S}^2$, where $\Omega^w_h$ \MMM is given by \EEE
$$\Omega_h^w:=w(\bar{\Omega}_h)\setminus w(\partial\bar{\Omega}_h)$$
which is well-defined, for \MMM $w$ \EEE is continuous. 
The magnetization $m$ is hence required to fulfill the {\it saturation}
constraint $|m|=1$ on $\Omega_h^w$.  \EEE

In what follows, \RRR for every $x\in\R^3$ in the referential space we write
$x=(x',x_3)$ \RRR where $x'\in\R^2$ is referred to as the planar coordinates of $x$, and we denote by $\nabla'$  the gradient with respect to such \EEE planar
coordinates. \UUU We use \RRR instead \UUU the symbol $\xi \in \R^3$ to indicate
variable\MMM s \EEE in \MMM the \EEE actual space. \EEE

Following the approach in \cite{jam-kin2,kruzik.stefanelli.zeman,rybka-luskin} we consider the total energy $I_h$, defined as
\begin{align}
\label{eq:def-Ih}I_h(w,m)&:=\int_{\Omega_h}W(\nabla w(x),m\circ w(x))dx+\alpha\int_{\Omega_h^{w}}|\nabla m(\xi)|^2 d\xi+\int_{\Omega_h}|\nabla^2w(x)|^p\,dx\\
\nonumber &\quad+\int_{\Omega_h}\Phi(\nabla w(x))dx+\frac{\mu_0}{2}\int_{\R^3}|\nabla u_m(\xi)|^2 d\xi.
\end{align}
In the formula above, $W:\mthree\times \mathbb{S}^2\to[0,+\infty)$ is the \emph{elastic energy} density associated to the plate,  which is a continuous function  satisfying \EEE the following assumptions:
\begin{align}
&\label{eq:H1}\text{(Coercivity)}\quad \qquad\qquad  \exists c>0 \text{
  such that }W(F,\lambda)\geq c|F|^p-\frac1c, \\
&\label{eq:H2}\text{(Frame indifference)}\quad \ W(RF,R\lambda)=W(F,\lambda),\\
&\label{eq:H3}\text{(Magnetic parity)}\quad \quad \,W(F,\lambda)=W(F,-\lambda)
%%&\label{eq:H3}\text{(Polyonvexity)}\quad \qquad \,\,W(F,\lambda)=\hat{W}(F,\textrm{cof}\,F, \lambda)
\end{align}
for every $F\in \mthree$,     $R\in {\rm SO}(3)$, \EEE and $\lambda\in \mathbb{S}^2$. % where $\hat{W}:\mthree\times \mathbb{S}^2\to[0,+\infty)$ is a continuous function, and is such that $\hat{W}(\cdot,\cdot,\lambda)$ is convex for every $\lambda\in \mathbb{S}^2$ (see \cite{ball}).
\UUU In fact, assumptions \eqref{eq:H2}-\eqref{eq:H3} are not strictly
needed for the analysis, but rather required by modeling considerations. \EEE

The second term in \UUU the expression of $I_h$ in \EEE
\eqref{eq:def-Ih} is the \emph{exchange energy}. The constant $\alpha$
is related to the size of ferromagnetic texture. \UUU The material is
assumed to be of {\it nonsimple} type \cite{Kruzik-Roubicek}. This is
expressed by the occurrence of the third term in $I_h$, providing a
higher-order contribution and a further length scale to the problem. \EEE  Regarding the fourth term, we will require that  $\Phi:\mthree\to [0,+\infty)$ is a continuous map satisfying the following assumptions
\begin{align}
&\nn \Phi(F)\to +\infty\quad\text{as }\,{\rm det}\,F\to 0^+,\\
&\nn\Phi(F)=+\infty\quad\text{if }\,{\rm det}\,F\leq 0,\\
&\Phi(F)\geq \frac{1}{C}({\rm det}\,F)^{-q}\quad\text{ for some $C>0$ \EEE and for every }\,F\in \mthree\,\text{with }{\rm det}\,F>0,\label{eq:hp-Phi}
\end{align}
where $q>\frac{3p}{p-3}$. \UUU This last quantification is \MMM introduced \EEE in
\EEE \cite{healey.kromer} and ensures that, \UUU for all \MMM $\lambda>0$ \EEE  and  $w\in W^{2,p}(\O_h;\R^3)$ such that  
\begin{align}\label{HK}
\int_{\Omega_h}|\nabla^2 w(x)|^p\,dx+\int_{\Omega_h}\Phi(\nabla w(x))dx < \MMM \lambda \EEE
\end{align} 
there exists $c>0$ depending on \MMM $\lambda>0$ \EEE with the property that 
$$%\be{eq:pos-det}
{\rm det }\,\nabla w>c\quad\text{ in }\bar\Omega_h.
$$%\ee 
Note that the left-hand side of \UUU inequality \EEE \eqref{HK} is a part of the energy functional \eqref{eq:def-Ih}.
\EEE
The last term in \eqref{eq:def-Ih} represents the \emph{magnetostatic energy}. In particular, $\mu_0$ is the \RRR \emph{permittivity} \EEE of void, and $u_m$ solves the Maxwell equation
$$%\be{eq:maxwell}
\nabla\cdot(-\mu_0\nabla u_m+\chi_{\Omega_h^w}m)=0\quad\text{in }\R^3,
$$%\ee
where $\chi_{\Omega_h^w}$ is the characteristic function of the set
$\Omega_h^w$. \UUU For simplicity, we assume \MMM that \EEE the deformations $\UUU w
\EEE $ satisfy the boundary conditions 
%\begin{equation}
%\label{eq:bc}
$$
\UUU w
\EEE =\textrm{id}\qquad\text{and}\qquad \nabla \UUU w
\EEE =\textrm{Id}\qquad\text{on }\,\partial\omega\times\big(-\tfrac12,\tfrac12\big).
$$%\end{equation}
\UUU To consider alternative boundary conditions would call for
solving some additional technicalities which, we believe, would
excessively complicate the argument. We hence leave this extension to
some possible further investigation. \EEE

\subsection{Change of variables}
As customary in dimension reduction, we perform the change of variables
$$\phi_h:\Omega\to\Omega_h,\quad\phi_h(x):=(x_1,x_2,hx_3)\quad\text{for a.e. }x\in\Omega.$$
Setting $y:=w\circ \phi_h$, \UUU $\Omega^y := y (\bar \Omega)\setminus
y(\partial \bar \Omega)$, \EEE and $E^h(y,m):=\frac{1}{h}I^h(w,m)$, we obtain
\begin{align*}
%\label{eq:def-eh}
E_h(y,m)&:=\int_{\Omega}W(\nh y(x),m\circ y(x))dx+\frac{\alpha}{h}\int_{\Omega^y}|\nabla m(\xi)|^2d\xi+\int_{\Omega}|\nabla^2_h y(x)|^p\,dx\\
\nonumber &\quad+\int_{\Omega}\Phi(\nabla_h y(x))\,dx+
\frac{\mu_0}{2h}\int_{\R^3}|\nabla u_m(\xi)|^2d\xi,
\end{align*}
where
$$\nabla\cdot(-\mu_0\nabla u_m+\chi_{\Omega^y}m)=0\quad\text{in
}\R^3.$$
\UUU Above, \EEE $\nh$ and $\nh^2$ are the differential operators defined as
$$\nh v :=\Big(\partial_1v \big|\partial_2 v \big|\frac{\partial_3 v }{h}\Big),\quad\text{and}\quad \nh^2 v:=\Bigg(\begin{array}{ccc}\qquad\partial^2_{11}v&\qquad\partial^2_{21}v&h^{-1}\partial^2_{31}v\\\qquad\partial^2_{12}v&\qquad\partial^2_{22}v&h^{-1}\partial^2_{32}v\\h^{-1}\partial^2_{13}v&h^{-1}\partial^2_{23}v&h^{-2}\partial^2_{33}v\end{array}\Bigg)$$
for every $v\in W^{2,p}(\Omega)$.
Note that the three-dimensional Ciarlet-Ne\v{c}as condition becomes
\be{eq:cn-rescaled}
\int_{\Omega}{\rm det}\,\nabla_h y\,dx\leq \frac{\mathcal{L}^3(\Omega^{y})}{h}.
\ee
Condition \eqref{eq:cn-rescaled}  \UUU provides scant information \EEE
in the thin-film regime, \UUU for \EEE
it leads to  the inequality 
\begin{align*}%\label{cn2D}
\int_\omega (\partial_1 y\times\partial_2 y)\cdot b\, {\rm d} S\le
\lim_{h\to 0} \frac{\mathcal{L}^3(\Omega^{y})}{h} 
\end{align*}
where $b$ is  a {\it Cosserat} vector obtained as  $b=\lim_{h\to 0} h^{-1} \partial_3 y^h$ in $ W^{1,p}(\omega;\R^3)$.
In particular, if $b= (\partial_1 y\times\partial_2 y)/|\partial_1 y\times\partial_2 y|$, i.e., it is the unit normal vector to the film in the deformed configuration,
and if $\lim_{h\to 0} \frac{\mathcal{L}^3(\Omega^{y})}{h}=\mathcal{H}^2(y(\omega))$ we get
\begin{align}\label{cn2D-spec}
\int_\omega  |\partial_1 y\times\partial_2 y|\,{\rm d}S\le \mathcal{H}^2(y(\omega)).
\end{align}
The left-hand side of \eqref{cn2D-spec} is the area of the deformed film calculated by the change-of-variables formula while the right-hand side is 
the measured area. Hence, \eqref{cn2D-spec} is violated by a {\it
  folding} deformation, which \UUU should be \EEE admissible \RRR among the family of realistic thin-film deformations,  while  \eqref{cn2D-spec} \EEE
is satisfied if the film crosses itself,  which  \UUU violates
non-self-interpenetration \RRR of matter \UUU and is hence not admissible. \EEE
On the other hand, if $y:\Omega\to\R^3$ is injective then  \eqref{cn2D-spec} is satisfied. The situation is depicted in the following figures\RRR, i.e., Figures \ref{fig:inj}-\ref{fig:non-inj}. 
%Figures \ref{fig:inj}, \ref{\ref,fig:carp}, and \ref{fig:non-inj}. 
\EEE
%\begin{comment}
\begin{center}
%flying carpet
\begin{figure}[h]
\begin{tikzpicture}
\begin{axis}[hide axis]
%%[
    %%title={$x \exp(-x^2-y^2)$}, 
    %%xlabel=$x$, ylabel=$y$,
	%%small,
%%]
\addplot3[
	surf,
	domain=-2:2,
	domain y=-1.3:1.3,
] 
	{exp(-x^2-y^4)*x};
\end{axis}
\end{tikzpicture}
\caption{An injective deformation satisfying  \eqref{cn2D-spec}}
\label{fig:inj} 
\end{figure}
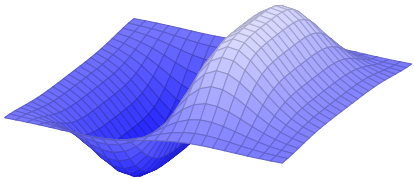
\hspace{2 cm}
%rolled-up carpet
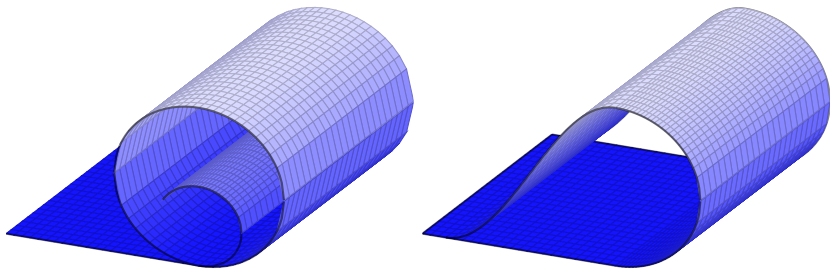
\begin{figure}[h]
\begin{tikzpicture}
\begin{axis}[domain=-(2):3.14, hide axis]
%\addplot3[surf] {min(0.,(1-0.3)*e^(-x*(y/100)*(1-0.3))-e^(-x*(y/100))};
\addplot3[surf, domain=-(2):1.57, domain y=-(1.57):1.57] {max(0.,0.)};
%\addplot3[domain=-(1.565):1.57,samples=80,samples y=0, mark=none,black, opacity=0.5,thick]({x},{1.57}, {max(0.,1.57+1.57*sin(deg(x)))});
\addplot3[domain=-(2):1.57, domain y=-(1.57):1.57,samples=80,samples y=0, mark=none,black, opacity=0.5,thick]({x},{-1.57}, {max(0.,0.)});
%\addplot3[domain=1.57:3.14, domain y=-(1.57):1.57,samples=80,samples y=0, mark=none,black, opacity=0.5,thick]({x},{1.57}, {max(1.57+sqrt((1.57)^2-(x-1.57)^2),1.57-sqrt((1.57)^2-(x-1.57)^2) )});
%\addplot3[domain=1.57:3.14, domain y=-(1.57):1.57,samples=80,samples y=0, mark=none,black, opacity=0.5,thick]({x},{1.57}, {min(1.57+sqrt((1.57)^2-(x-1.57)^2),1.57-sqrt((1.57)^2-(x-1.57)^2) )});

%\addplot3[surf, domain=-(1.57):1.57,] {max(0.,1.57+1.57*sin(deg(x)))};
\addplot3[surf, domain=1.57:2.355, domain y=-(1.57):1.57,] {min(0.785+sqrt((0.785)^2-(x-1.57)^2),0.785-sqrt((0.785)^2-(x-1.57)^2) )};

\addplot3[surf, domain=0.:3.14, domain y=-(1.57):1.57,] {min(1.57+sqrt((1.57)^2-(x-1.57)^2),1.57-sqrt((1.57)^2-(x-1.57)^2) )};
\addplot3[surf, domain=0.9:2.355, domain y=-(1.57):1.57,] {max(0.785+sqrt((0.785)^2-(x-1.57)^2),0.785-sqrt((0.785)^2-(x-1.57)^2) )};
\addplot3[surf, domain=0.:3.14, domain y=-(1.57):1.57,] {max(1.57+sqrt((1.57)^2-(x-1.57)^2),1.57-sqrt((1.57)^2-(x-1.57)^2) )};
\addplot3[domain=1.57:2.355, domain y=-(1.57):1.57,samples=80,samples y=0, mark=none,black, opacity=0.5,thick] ({x},{-(1.57)},{min(0.785+sqrt((0.785)^2-(x-1.57)^2),0.785-sqrt((0.785)^2-(x-1.57)^2) )});
\addplot3[domain=0.9:2.355, domain y=-(1.57):1.57,samples=80,samples y=0, mark=none,black, opacity=0.5,thick] ({x},{-(1.57)},{max(0.785+sqrt((0.785)^2-(x-1.57)^2),0.785-sqrt((0.785)^2-(x-1.57)^2) )});
\addplot3[domain=0.:3.14, domain y=-(1.57):1.57,samples=80,samples y=0, mark=none,black, opacity=0.5,thick] ({x},{-(1.57)},{min(1.57+sqrt((1.57)^2-(x-1.57)^2),1.57-sqrt((1.57)^2-(x-1.57)^2) )});
\addplot3[domain=0.:3.14, domain y=-(1.57):1.57,samples=80,samples y=0, mark=none,black, opacity=0.5,thick] ({x},{-(1.57)},{max(1.57+sqrt((1.57)^2-(x-1.57)^2),1.57-sqrt((1.57)^2-(x-1.57)^2) )});
%\addplot3[domain=-(1.565):1.57,samples=80,samples y=0, mark=none,black, opacity=0.5,thick]({x},{-1.57}, {max(0.,1.57+1.57*sin(deg(x)))});
%\addplot3[domain=1.57:3.14, domain y=-(1.57):1.57,samples=80,samples y=0, mark=none,black, opacity=0.5,thick]({x},{-1.57}, {max(1.57+sqrt((1.57)^2-(x-1.57)^2),1.57-sqrt((1.57)^2-(x-1.57)^2) )});
%\addplot3[domain=1.57:3.14, domain y=-(1.57):1.57,samples=80,samples y=0, mark=none,black, opacity=0.5,thick]({x},{-1.57}, {min(1.57+sqrt((1.57)^2-(x-1.57)^2),1.57-sqrt((1.57)^2-(x-1.57)^2) )});
%\addplot3[domain=-(2):1.57, domain y=-(1.57):1.57,samples=80,samples y=0, mark=none,black, opacity=0.5,thick]({x},{-1.57}, {max(0.,0.)});
%\addplot3[domain=-(2):1.57, domain y=-(1.57):1.57,samples=80, mark=none,black, opacity=0.5,thick]({-2},{y}, {max(0.,0.)});

%\addplot3[domain=-(1.57):1.57, samples=80, mark=none, black, opacity=0.5,thick]({x},{y},{max(0.,1.57+1.57*sin(deg(x)))});
\end{axis}

 \end{tikzpicture}
\begin{tikzpicture}
\begin{axis}[domain=-(2):3.14, hide axis]
%\addplot3[surf] {min(0.,(1-0.3)*e^(-x*(y/100)*(1-0.3))-e^(-x*(y/100))};
\addplot3[surf, domain=-(2):1.57, domain y=-(1.57):1.57] {max(0.,0.)};
\addplot3[domain=-(1.565):1.57,samples=80,samples y=0, mark=none,black, opacity=0.5,thick]({x},{1.57}, {max(0.,1.57+1.57*sin(deg(x)))});
\addplot3[domain=-(2):1.57, domain y=-(1.57):1.57,samples=80,samples y=0, mark=none,black, opacity=0.5,thick]({x},{1.57}, {max(0.,0.)});
\addplot3[domain=1.57:3.14, domain y=-(1.57):1.57,samples=80,samples y=0, mark=none,black, opacity=0.5,thick]({x},{1.57}, {max(1.57+sqrt((1.57)^2-(x-1.57)^2),1.57-sqrt((1.57)^2-(x-1.57)^2) )});
\addplot3[domain=1.57:3.14, domain y=-(1.57):1.57,samples=80,samples y=0, mark=none,black, opacity=0.5,thick]({x},{1.57}, {min(1.57+sqrt((1.57)^2-(x-1.57)^2),1.57-sqrt((1.57)^2-(x-1.57)^2) )});

\addplot3[surf, domain=-(1.57):1.57,] {max(0.,1.57+1.57*sin(deg(x)))};
\addplot3[surf, domain=1.57:3.14, domain y=-(1.57):1.57,] {min(1.57+sqrt((1.57)^2-(x-1.57)^2),1.57-sqrt((1.57)^2-(x-1.57)^2) )};
\addplot3[surf, domain=1.57:3.14, domain y=-(1.57):1.57,] {max(1.57+sqrt((1.57)^2-(x-1.57)^2),1.57-sqrt((1.57)^2-(x-1.57)^2) )};
\addplot3[domain=-(1.565):1.57,samples=80,samples y=0, mark=none,black, opacity=0.5,thick]({x},{-1.57}, {max(0.,1.57+1.57*sin(deg(x)))});
\addplot3[domain=1.57:3.14, domain y=-(1.57):1.57,samples=80,samples y=0, mark=none,black, opacity=0.5,thick]({x},{-1.57}, {max(1.57+sqrt((1.57)^2-(x-1.57)^2),1.57-sqrt((1.57)^2-(x-1.57)^2) )});
\addplot3[domain=1.57:3.14, domain y=-(1.57):1.57,samples=80,samples y=0, mark=none,black, opacity=0.5,thick]({x},{-1.57}, {min(1.57+sqrt((1.57)^2-(x-1.57)^2),1.57-sqrt((1.57)^2-(x-1.57)^2) )});
\addplot3[domain=-(2):1.57, domain y=-(1.57):1.57,samples=80,samples y=0, mark=none,black, opacity=0.5,thick]({x},{-1.57}, {max(0.,0.)});
\addplot3[domain=-(2):1.57, domain y=-(1.57):1.57,samples=80, mark=none,black, opacity=0.5,thick]({-2},{y}, {max(0.,0.)});

%\addplot3[domain=-(1.57):1.57, samples=80, mark=none, black, opacity=0.5,thick]({x},{y},{max(0.,1.57+1.57*sin(deg(x)))});
\end{axis}
%\begin{axis}[domain=3.14:4.71, hide axis]
%
%\addplot3[domain=4:30,samples=80,samples y=0,mark=none,black, opacity=0.5,thick]({x},{118.89/x},{0.});
%\addplot3[domain=0:30,samples=80,samples y=0,mark=none,black, opacity=0.5,thick]({x},{30.},{max(0.,(1-0.3)*e^(-x*(30./100)*(1-0.3))-e^(-x*(30./100)))});
%\addplot3[domain=0:30,samples=80,samples y=0,mark=none,black, opacity=0.5,thick]({x},{0.},{max(0.,(1-0.3)*e^(-x*(0./100)*(1-0.3))-e^(-x*(0./100)))});
%\addplot3[domain=0:30,samples=80,samples y=0,mark=none,black, opacity=0.5,thick]({x},{0.},{min(0.,(1-0.3)*e^(-x*(0./100)*(1-0.3))-e^(-x*(0./100)))});
%\addplot3[domain=0:30,samples=80,samples y=0,mark=none,black, opacity=0.5,thick]({0.},{x},{max(0.,(1-0.3)*e^(-0.*(x/100)*(1-0.3))-e^(-0.*(x/100)))});
%\addplot3[domain=0:30,samples=80,samples y=0,mark=none,black, opacity=0.5,thick]({30.},{x},{max(0.,(1-0.3)*e^(-30.*(x/100)*(1-0.3))-e^(-30.*(x/100)))});
%\addplot3[domain=0:30,samples=80,samples y=0,mark=none,white, opacity=0.5,thick]({30.},{x},{min(0.,(1-0.3)*e^(-30.*(x/100)*(1-0.3))-e^(-30.*(x/100)))});
%\end{axis}

 \end{tikzpicture}
\caption{Two deformations not satisfying  \eqref{cn2D-spec} in the regions in which a self-contact occurs}
\label{fig:carp}
\end{figure}

%self intersection
\begin{figure}[h]
\begin{tikzpicture}
\begin{axis}[domain=-(2):3.14, hide axis]
%\addplot3[surf] {min(0.,(1-0.3)*e^(-x*(y/100)*(1-0.3))-e^(-x*(y/100))};
\addplot3[surf, domain=0:3.14, domain y=-(1.57):1.57,] {min(1.57+sqrt((3.14)^2-(x-3.14)^2),1.57-sqrt((3.14)^2-(x-3.14)^2) )};
\addplot3[domain=0:3.14, domain y=-(1.57):1.57,samples=80,samples y=0, mark=none,black, opacity=0.5,thick] ({x},{-1.57},{min(1.57+sqrt((3.14)^2-(x-3.14)^2),1.57-sqrt((3.14)^2-(x-3.14)^2) )});
\addplot3[domain=0:3.14, domain y=-(1.57):1.57,samples=80,samples y=0, mark=none,black, opacity=0.5,thick] ({x},{1.57},{min(1.57+sqrt((3.14)^2-(x-3.14)^2),1.57-sqrt((3.14)^2-(x-3.14)^2) )});

\addplot3[surf, domain=-(2):1.57, domain y=-(1.57):1.57] {max(0.,0.)};
\addplot3[domain=-(2):1.57, domain y=-(1.57):1.57,samples=80,samples y=0, mark=none,black, opacity=0.5,thick] ({x},{1.57},{0.});
\addplot3[domain=0:0.6, domain y=-(1.57):1.57,samples=80,samples y=0, mark=none,black, thick, dashed] ({x},{-1.57},{min(1.57+sqrt((3.14)^2-(x-3.14)^2),1.57-sqrt((3.14)^2-(x-3.14)^2) )});

\addplot3[domain=-(2):1.57, domain y=-(1.57):1.57,samples=80,samples y=0, mark=none,black, opacity=0.5,thick]({x},{-1.57}, {max(0.,0.)});
\addplot3[domain=0:3.14, domain y=-(1.57):1.57,samples=80,mark=none,black, very thick, dashed] ({0},{y},{1.57});
\addplot3[domain=1.57:3.14, domain y=-(1.57):1.57,samples=80,samples y=0, mark=none,black, opacity=0.5,thick]({x},{1.57}, {max(1.57+sqrt((1.57)^2-(x-1.57)^2),1.57-sqrt((1.57)^2-(x-1.57)^2) )});
\addplot3[domain=1.57:3.14, domain y=-(1.57):1.57,samples=80,samples y=0, mark=none,black, opacity=0.5,thick]({x},{1.57}, {min(1.57+sqrt((1.57)^2-(x-1.57)^2),1.57-sqrt((1.57)^2-(x-1.57)^2) )});

\addplot3[surf, domain=1.57:3.14, domain y=-(1.57):1.57,] {min(1.57+sqrt((1.57)^2-(x-1.57)^2),1.57-sqrt((1.57)^2-(x-1.57)^2) )};

\addplot3[surf, domain=0.:3.14, domain y=-(1.57):1.57,] {max(1.57+sqrt((1.57)^2-(x-1.57)^2),1.57-sqrt((1.57)^2-(x-1.57)^2) )};

%\addplot3[domain=1.57:2.355, domain y=-(1.57):1.57,samples=80,samples y=0, mark=none,black, opacity=0.5,thick] ({x},{-(1.57)},{min(0.785+sqrt((0.785)^2-(x-1.57)^2),0.785-sqrt((0.785)^2-(x-1.57)^2) )});
%\addplot3[domain=0.9:2.355, domain y=-(1.57):1.57,samples=80,samples y=0, mark=none,black, opacity=0.5,thick] ({x},{-(1.57)},{max(0.785+sqrt((0.785)^2-(x-1.57)^2),0.785-sqrt((0.785)^2-(x-1.57)^2) )});
\addplot3[domain=1.57:3.14, domain y=-(1.57):1.57,samples=80,samples y=0, mark=none,black, opacity=0.5,thick] ({x},{-(1.57)},{min(1.57+sqrt((1.57)^2-(x-1.57)^2),1.57-sqrt((1.57)^2-(x-1.57)^2) )});
\addplot3[domain=0.:3.14, domain y=-(1.57):1.57,samples=80,samples y=0, mark=none,black, opacity=0.5,thick] ({x},{-(1.57)},{max(1.57+sqrt((1.57)^2-(x-1.57)^2),1.57-sqrt((1.57)^2-(x-1.57)^2) )});

%\addplot3[domain=1.57:3.14, domain y=-(1.57):1.57,samples=80,samples y=0, mark=none,black, opacity=0.5,thick]({x},{-1.57}, {max(1.57+sqrt((1.57)^2-(x-1.57)^2),1.57-sqrt((1.57)^2-(x-1.57)^2) )});
%\addplot3[domain=1.57:3.14, domain y=-(1.57):1.57,samples=80,samples y=0, mark=none,black, opacity=0.5,thick]({x},{-1.57}, {min(1.57+sqrt((1.57)^2-(x-1.57)^2),1.57-sqrt((1.57)^2-(x-1.57)^2) )});
%\addplot3[domain=-(2):1.57, domain y=-(1.57):1.57,samples=80,samples y=0, mark=none,black, opacity=0.5,thick]({x},{-1.57}, {max(0.,0.)});
%\addplot3[domain=-(2):1.57, domain y=-(1.57):1.57,samples=80, mark=none,black, opacity=0.5,thick]({-2},{y}, {max(0.,0.)});

%\addplot3[domain=-(1.57):1.57, samples=80, mark=none, black, opacity=0.5,thick]({x},{y},{max(0.,1.57+1.57*sin(deg(x)))});
\end{axis}

 \end{tikzpicture}
 \caption{A \UUU self-interpenetrating \EEE deformation satisfying  \eqref{cn2D-spec}} 
 \label{fig:non-inj}
 \end{figure}
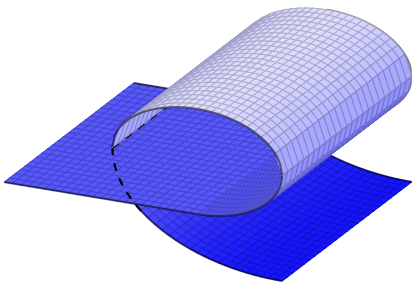
\end{center}

%\end{comment}

\EEE

%We point out that if $y$ is such that $E_h(y,m)<+\infty$, then the set $\Omega^y$ has the cone property. Indeed, by the definition of $E_h$ and the boundedness of the terms involving the gradient of the cofactor and the determinant of $\nabla_h y$ it follows that $y\in W^{1,\infty}(\Omega;\R^3)$, and by \eqref{eq:cn-rescaled} we get that it is also globally injective almost everywhere. 

In what follows we analyze the asymptotic behaviour of sequences $(y^h,m^h)\in W^{2,p}(\Omega;\R^3)\times W^{1,2}(\Omega^{y^h};\R^3)$ satisfying the uniform energy estimate
\be{eq:uniform-en-estimate}
E_h(y^h,m^h)\leq C,
\ee
and the boundary conditions 
\begin{equation}
\label{eq:rescaled-bc}
y(x)=(x',hx_3)\qquad\text{and}\qquad \nabla y(x)=\Bigg(\begin{array}{ccc}1&0&0\\0&1&0\\0&0&h\end{array}\Bigg)\qquad\text{on }\,\partial\omega\times\big(-\tfrac h2,\tfrac h2\big).
\end{equation}
\UUU A {\it caveat} on notation: in \eqref{eq:uniform-en-estimate} and
in the following the symbol $C$ is used to denote a generic constant \RRR that may possibly change from line to line and that  always depends only on model data and not on $h$. \EEE
%depending on data but independent of $h$, which may possibly change from line to line. \EEE

We point out that, without the \UUU $\Phi$ term in the energy \EEE and
the Ciarlet-Ne\v{c}as condition, constants deformations $y$ having
null gradient, null hessian, and such that the measure of the deformed
set is zero (so that the exchange energy gives no contribution) would
be energetically favourable both for the elastic and the exchange
energy. The associated magnetic field would then concentrate in a
point. The  \UUU $\Phi$ term \EEE  in our model prevent this degenerate situation from happening.

%%%%%%%%%%%%%%%%%%%%%%%%%%%%%%%%%%%%%%%%%%%%
\section{\UUU Main results \EEE}%$\Gamma$-convergence}
\label{sec:gamma}
\UUU This section is devoted to the statement of our main
$\Gamma$-convergence results. All proofs are postponed to the
following Section \ref{sec:proofs}. \EEE

 For notational convenience, for every open set $U\subset \R^2$ we
 denote by $\mathring{W}^{k,p}(U;\R^n)$ the set of $W^{k,p}$-maps
 having \UUU zero mean \EEE on $U$, i.e., $y\in \mathring{W}^{k,p}(U;\R^n)$ if $y\in W^{k,p}(U;\R^n)$ and  $\int_U y(x')\,\md x'=0$. As it is standard, we write $\mathring{W}^{k,p}(U)$ if $n=1$. \EEE

We first introduce the \RRR set  $\mathcal{A}$ \EEE of \emph{admissible
  limiting deformations  $y:\omega\to\R^3$, Cosserat vectors
  $b:\omega\to\R^3$, and magnetizations $\mathscr{M}:\UUU \omega \to
  \mathbb S^2$}, \EEE defined as
\begin{align*}
\mathcal{A}:=\{&(y,b,\mathscr{M}):\,y\in \mathring{W}^{2,p}(\omega;\R^3),\,b\in W^{1,p}(\omega;\R^3),\mathscr{M}\in W^{1,2}(\omega;\mathbb{S}^2),\\
&%\label{eq:bc-liminf} 
y=\text{id}\text{ and }\big(\nabla' y|b\big)=\text{Id}\text{ on }\,\partial\omega,\\
&%\label{eq:limit-cont} 
(\nabla'y|b)^{-1}\in
C^0(\bar{\omega};\mthree),\ \text{det}(\nabla'y|b)\in
C^0(\bar{\Omega}), \ \text{and}\ \text{det}(\nabla'y|b)>\ep  \  \text{for some }\ep>0\}.
\end{align*}

\UUU Let us \MMM first \UUU state the following lemma, which \EEE will be
instrumental in characterizing the limiting stray fields \UUU and
formulating the limiting functional. \UUU As mentioned, the lemma is proved in Section \ref{sec:proofs}
below. \EEE

\begin{lemma}
\label{lem:stray}
Let $(y,b,\mathscr{M})\in \mathcal{A}$. Denote by $\widetilde{(\nabla'y|b)}$ and $\bar{\mathscr{M}}$ the quantities
 \begin{equation}
 \label{eq:def-grad-ex}
 \widetilde{(\nabla'y|b)}(x'):=
 \begin{cases}
 (\nabla'y|b)(x')&\text{if}\ x'\in \omega\\
 {\rm Id}&\UUU \text{if}\ x'\in \mathbb{R}^2\setminus\omega. \EEE%\text{otherwise in}\,\mathbb{R}^2.
 \end{cases}
 \end{equation}
 and
 \begin{equation}
 \label{eq:def-m-bar}
 \bar{\mathscr{M}}(x'):=
 \begin{cases}
 \mathscr{M}(x')&\text{if}\ x'\in\omega\\
 0&\UUU \text{if}\ x'\in \mathbb{R}^2\setminus\omega. \EEE%\text{otherwise in}\,\mathbb{R}^2.
 \end{cases}
 \end{equation}
Then, the system
\begin{align}
&\label{eq:M1}\Big\{{\rm cof}\,\widetilde{(\nabla' y|b)^\top}\Big[\mu_0 \widetilde{(\nabla'y|b)}^{-T}\Big(\nabla'\mathscr{U}|\mathscr{V}\Big)^T-\bar{\mathscr{M}}\Big]\Big\}_3=0\quad\text{in }\,\R^2,\\
&\label{eq:M2} {\rm div}_{x'}\,\left\{\begin{array}{c}\Big[{\rm cof}\,\widetilde{(\nabla' y|b)^\top}\Big(\mu_0 \widetilde{(\nabla'y|b)}^{-T}\Big(\nabla'\mathscr{U}|\mathscr{V}\Big)^T-\bar{\mathscr{M}}\Big)\Big]_1\\
\Big[{\rm cof}\,\widetilde{(\nabla' y|b)^\top}\Big(\mu_0 \widetilde{(\nabla'y|b)}^{-T}\Big(\nabla'\mathscr{U}|\mathscr{V}\Big)^T-\bar{\mathscr{M}}\Big)\Big]_2 \end{array}\right\}=0\quad\text{in }\,\R^2,
\end{align}
has a unique solution $(\mathscr{U},\mathscr{V})\in W^{1,2}(\R^2)\times L^2(\R^2)$ satisfying $\int_{\omega}\mathscr{U}\,\md x'=0$.
\end{lemma}

\UUU The \EEE limiting energy is \UUU given \EEE  by the functional
$\UUU \mathcal{F} \EEE :\mathcal{A}\to [0,+\infty)$ \RRR defined as \EEE
\begin{align*}
\nn\mathcal{F}(y,b,\mathcal{M})&:=\int_{\omega}W\big((\nabla' y|b),\,\mathscr{M}\big)\,\md x'+ \alpha\int_{\Omega}|(\nabla'y|b)^{-T}(\nabla'\mathscr{M}|0)|^2{\rm det}\,(\nabla'y|b)\,\md x\\
&\nn\quad+\int_{\omega}(|(\nabla')^2 y|^2+2|\nabla' b|^2)^{p/2}\,\md x'+\int_{\omega}\Phi(\nabla'y|b)\,\md x'\\
&%\label{eq:lim-en}
\quad+\frac{\mu_0}{2}\int_{\omega}{\rm cof}\,(\nabla'y|b)\mathscr{M}\cdot \Big(\nabla'\mathscr{U}_{y,b,\mathscr{M}}|\mathscr{V}_{y,b,\mathscr{M}}\Big)^T\,\md x'
\end{align*}
 for every $(y,b,\mathscr{M})\in \mathcal{A}$, where the pair $\Big(\mathscr{U}_{y,b,\mathscr{M}}|\mathscr{V}_{y,b,\mathscr{M}}\Big)\in \mathring{W}^{1,2}(\omega)\times L^2(\omega)$ is the restriction to $\omega$ of the unique solution to \eqref{eq:M1}--\eqref{eq:M2} in the sense of Lemma \ref{lem:stray}.

We start by providing a lower bound for the asymptotic behavior of the
functionals $\{E_h\}_h$ along sequences of deformations and
magnetizations with equibounded energies. \UUU Again, \RRR the \UUU proof is
postponed to Section \ref{sec:proofs}. \EEE
\begin{theorem}[Compactness and \UUU $\Gamma$-$\liminf$ \EEE inequality]
\label{prop:first-liminf}
Let $\{(y^h,m^h)\}\subset W^{2,p}(\Omega;\R^3)\times
W^{1,2}(\Omega^{y^h};\R^3)$ be such that
\eqref{eq:uniform-en-estimate} holds true. Then, there exist
$(y,b,\mathscr{M})\in \mathcal{A}$ \UUU and \EEE $d\in
L^p(\Omega;\UUU\R^3 \EEE)$ such that up to the extraction of a (not relabeled) subsequence there holds
\begin{align}
&\label{eq:wk-null-av-def} y^h\wk y\quad\text{weakly in }W^{2,p}(\Omega;\R^3),\\
&\label{eq:wk-conv-grad} \nabla_h y^h\wk (\nabla' y|b)\quad\text{weakly in }W^{1,p}(\Omega;\mthree),\\
&\label{eq:weak-d33} \frac{\partial^2_{33}y^h}{h^2}\wk d\quad\text{weakly in }L^p(\Omega;\R^3).
\end{align}
Additionally, there exist $\eta\in L^2(\Omega;\R^3)$ and $\mathscr{V}\in L^2(\Omega)$ such that $\int_{-\tfrac12}^{\tfrac12}\mathscr{V}\,\md x_3=\mathscr{V}_{y,b,\mathscr{M}}$, and up to subsequences we have
\begin{align}
&\label{eq:m-wk} m^h\circ y^h\wk\mathscr{M}\quad\text{weakly in }W^{1,2}(\Omega;\R^3),\\
&\label{eq:m-wk-bis} \nabla_h(m^h\circ y^h)\wk (\nabla'\mathscr{M}|\eta)\quad\text{weakly in }L^{2}(\Omega;\mthree),\\
&\label{eq:vh-one} u_{m^h}\circ y^h-\fint_{\Omega}u_{m^h}\circ y^h\,\md x\wk\mathscr{U}_{y,b,\mathscr{M}}\quad\text{weakly in }W^{1,2}(\omega),\\
&\label{eq:vh-wk} \nabla_h (u_{m^h}\circ y^h)\wk (\nabla'\mathscr{U}_{y,b,\mathscr{M}}|\mathscr{V})^T\quad\text{weakly in }L^2(\Omega;\R^3).
\end{align}

Eventually, the following liminf inequality for the energy holds true:
\begin{align}
\label{eq:liminf-inequality}
\liminf_{h\to 0} E_h(y^h,m^h)&\geq \mathcal{F}(y,b,\mathscr{M}).
\end{align}
\end{theorem}

\UUU The statement of our second main result requires the
specification of the class of admissible deformations. This is given
through the following definition. \EEE

% Before stating \UUU a \EEE second main, we introduce the class $\mathcal{Y}$ of \emph{approximately injective deformations}. 

 \begin{definition}[Approximately injective deformations]\label{AID}
\UUU We define the set $\mathcal{Y}$ of  \emph{approximately injective
  deformations} as \EEE 
 \begin{align*}
 \mathcal{Y}:=\Big\{&y\in W^{2,p}(\omega;\R^3):\,\text{there exist }b\in W^{1,p}(\omega;\R^3)\text{ and }\mathscr{M}\in W^{1,2}(\omega;\mathbb{S}^2)\text{ such that }(y,b,\mathscr{M})\in \mathcal{A},\\
 &\text{and there exists a sequence }\{f_h\}_h\subset  W^{2,p}(\omega;\R^3)\text{ for which }\\
 &y^h(x):=y(x')+hx_3 b(x')+ f^h(x')\EEE \text{ satisfy }\eqref{eq:ciarlet-necas}\text{ and }h^{-2}f^h\to 0\text{ strongly in }W^{2,p}(\omega;\R^3)\text{ as }h\to 0\Big\}.
 \end{align*}
\end{definition}
  The deformations in Figure \ref{fig:inj} and on the right of Figure
  \ref{fig:carp} fulfill the requirements of Definition \ref{AID},
  whereas those depicted on the left of Figure \ref{fig:carp} and in
  Figure \ref{fig:non-inj} are not included in the above setting. \UUU
  Let us note that, \EEE although still not covering all realistic
  thin-film deformations, \UUU the set of approximately injective
  deformations \EEE encompasses a wider range of scenarios compared to those allowed by \eqref{cn2D-spec}. \EEE
 
 We provide below a construction of a recovery sequence for triples
 $(y,b,\mathscr{M})\in \mathcal{A}$ \UUU under the assumption that \EEE $y\in \mathcal{Y}$.

 \begin{theorem}[Optimality of the lower bound for approximately injective deformations]
 \label{prop:limsup}
Let \UUU $y\in  \mathcal{Y}$ and $b$ and $\mathscr{M}$ given by the
definition of $\mathcal{Y}$ so that $(y,b,\mathscr{M})\in
\mathcal{A}$. \EEE Then,
there exist\MMM s \EEE \UUU a recovery sequence \EEE $\{(y^h,m^h)\}_h\subset W^{2,p}(\Omega;\R^3)\times
W^{1,2}(\Omega^{y^h};\R^3)$ such that, setting $u^h$ as the solution
to the Maxwell \UUU system \EEE equation
$${\rm div}\,(-\mu_0 \nabla u_{m^h}+\chi_{\Omega^{y^h}}m^h)=0$$
\UUU with zero mean, \EEE there holds
\begin{align*}
&%\label{eq:wk-null-av-def-s}
 y^h-\fint_{\Omega}{y^h}\,\md x\to y\quad\text{strongly in }W^{2,p}(\Omega;\R^3),\\
&%\label{eq:wk-conv-grad-s} 
\nabla_h y^h\to (\nabla' y|b)\quad\text{strongly in }W^{1,p}(\Omega;\mthree),\\
&%\label{eq:weak-d33-s}
 \frac{\partial_{33}y^h}{h}\to 0\quad\text{strongly in }L^p(\Omega;\R^3).
\end{align*}
Additionally,
\begin{align*}
&%\label{eq:m-wk-s} 
m^h\circ y^h=\mathscr{M}\quad\text{for every }h>0,\\
&%\label{eq:m-wk-bis-s} 
\nabla_h(m^h\circ y^h)=(\nabla'\mathscr{M}|0)\quad\text{for every }h>0,\\
&%\label{eq:vh-one-s} 
u_{m^h}\circ y^h-\fint_{\Omega}u_{m^h}\circ y^h\,\md x\wk\mathscr{U}_{y,b,\mathscr{M}}\quad\text{weakly in }W^{1,2}(\omega;),\\
&%\label{eq:vh-wk-s} 
\nabla_h (u_{m^h}\circ y^h)\wk (\nabla'\mathscr{U}_{y,b,\mathscr{M}}|\mathscr{V}_{y,b,\mathscr{M}})^T\quad\text{weakly in }L^2(\Omega;\R^3),
\end{align*}
and the following limsup inequality for the energy holds true:
\begin{align*}
%\label{eq:limsup-inequality}
\limsup_{h\to 0} E_h(y^h,m^h)\leq \mathcal{F}(y,b,\mathscr{M}).
\end{align*}
\end{theorem}

\UUU In order to give a full $\Gamma$-convergence result, in the
remainder of the section we \EEE
% %%%%%%%%%%%%%%%%%%%%%%%%%%%%%%%%%%%%%%%%%%%%%%%%%%%%%%
% \section{The averaged invertibility constraint}
% \label{sec:add-gamma}
% In this section we
restrict our analysis to deformations satisfying the following \emph{uniform averaged invertibility constraint}: there exists a constant $C>0$ such that
\begin{equation}
\label{eq:average-inv}
\Big|\int_{-\frac12}^{\frac12} y^h(x',x_3)\,\md x_3-\int_{-\frac12}^{\frac12} y^h(z',x_3)\,\md x_3\Big|\geq C|x'-z'|\quad\text{for every}\,h>0,
\end{equation}
for every $x',\,z'\in \omega$. Note that the condition above has a
pointwise meaning because maps with uniformly bounded energies are
\UUU at least \EEE $C^1$-regular.

The key idea of \eqref{eq:average-inv} is that deformed vertical
fibers might intersect, but are, in average, distant enough, compared
to the distance of the original points in the cross section. 

\UUU Let
us start by remarking that, \EEE under the same assumptions of Proposition \ref{prop:first-liminf}, and assuming additionally \eqref{eq:average-inv}, the limiting deformations $y\in W^{2,p}(\omega;\R^3)$  have the additional property:
\begin{align}
&\label{eq:py} \text{There exists a constant $C>0$ such that $|y(x')-y(z')|\geq C|x'-z'|$ for every $x',\,z'\in \omega$}.
%&\label{eq:pM} \text{There ex
\end{align}
\UUU In fact, \EEE property \eqref{eq:py} follows from
\eqref{eq:wk-null-av-def} and \eqref{eq:average-inv}.
% \begin{proposition}[Characterization of the domain of the $\Gamma$-limit]
% \label{prop:lim-inj}
% Under the same assumptions of Proposition \ref{prop:first-liminf}, and assuming additionally \eqref{eq:average-inv}, the limiting deformations $y\in W^{2,p}(\omega;\R^3)$  have the additional property:
% \begin{align}
% &\label{eq:py} \text{There exists a constant $C>0$ such that $|y(x')-y(z')|\geq C|x'-z'|$ for every $x',\,z'\in \omega$}.
% %&\label{eq:pM} \text{There exists $m\in W^{1,2}(\omega;\mathbb{S}^2)$ such that $\mathscr{M}=m\circ y$ in $\omega$}.
% \end{align}
% \end{proposition}
% \begin{proof}
% Property \eqref{eq:py} follows from \eqref{eq:wk-null-av-def} and \eqref{eq:average-inv}. %Assertion \eqref{eq:pM} is a direct consequence of \eqref{eq:py}.
% \end{proof}
 In view of \UUU \eqref{eq:py} we are in the position of obtaining the
 following \EEE $\Gamma$-convergence result.
 
 \begin{theorem}[\UUU $\Gamma$-convergence under uniform averaged invertibility]
 \label{thm:complete-gamma-conv}
 Let $\{(y^h,m^h)\}\subset W^{2,p}(\Omega;\R^3)\times W^{1,2}(\Omega^{y^h};\R^3)$ be such that \eqref{eq:uniform-en-estimate} and \eqref{eq:average-inv} hold true. Then, there exist $(y,b,\mathscr{M})\in \mathcal{A}$, $d\in L^p(\Omega;\UUU\R^3\EEE)$, and $\ep>0$ satisfying \eqref{eq:py}, such that, up to the extraction of a (not relabeled) subsequence, there holds
\begin{align*}
&\nn y^h\wk y\quad\text{weakly in }W^{2,p}(\Omega;\R^3),\\
&\nn \nabla_h y^h\wk (\nabla' y|b)\quad\text{weakly in }W^{1,p}(\Omega;\mthree),\\
&\nn\frac{\partial^2_{33}y^h}{h^2}\wk d\quad\text{weakly in }L^p(\Omega;\R^3).
\end{align*}
Additionally, there exist $\eta\in L^2(\Omega;\R^3)$, and $\mathscr{V}\in L^2(\Omega)$ such that $\int_{-\tfrac12}^{\tfrac12}\mathscr{V}\,\md x_3=\mathscr{V}_{y,b,\mathscr{M}}$, and up to subsequences we have
\begin{align*}
&\nn m^h\circ y^h\wk\mathscr{M}\quad\text{weakly in }W^{1,2}(\Omega;\R^3),\\
&\nn \nabla_h(m^h\circ y^h)\wk (\nabla'\mathscr{M}|\eta)\quad\text{weakly in }L^{2}(\Omega;\mthree),\\
&\nn u_{m^h}\circ y^h-\fint_{\Omega}u_{m^h}\circ y^h\,\md x\wk\mathscr{U}_{y,b,\mathscr{M}}\quad\text{weakly in }W^{1,2}(\omega),\\
&\nn \nabla_h (u_{m^h}\circ y^h)\wk (\nabla'\mathscr{U}_{y,b,\mathscr{M}}|\mathscr{V})^T\quad\text{weakly in }L^2(\Omega;\R^3).
\end{align*}

Eventually, the following liminf inequality for the energy holds true:
\begin{align*}
\liminf_{h\to 0} E_h(y^h,m^h)&\geq \mathcal{F}(y,b,\mathscr{M}).
\end{align*}

Conversely, for every $(y,b,\mathscr{M})\in \mathcal{A}$ with $y$ satisfying \eqref{eq:py} there exist $\{(\bar{y}^h,\bar{m}^h)\}_h\subset W^{2,p}(\Omega;\R^3)\times W^{1,2}(\Omega^{y^h};\R^3)$ such that, setting $u_{\bar{m}^h}$ as the solution to the Maxwell's equation
$${\rm div}\,(-\mu_0 \nabla u_{\bar{m}^h}+\chi_{\Omega^{\bar{y}^h}}\bar{m}^h)=0$$
\UUU with zero mean, \EEE there holds
\begin{align*}
& \bar{y}^h-\fint_{\Omega}{\bar{y}^h}\,\md x\to y\quad\text{strongly in }W^{2,p}(\Omega;\R^3),\\
&\nabla_h \bar{y}^h\to (\nabla' y|b)\quad\text{strongly in }W^{1,p}(\Omega;\mthree),\\
& \frac{\partial_{33}\bar{y}^h}{h}\to 0\quad\text{strongly in }L^p(\Omega;\R^3).
\end{align*}
Additionally,
\begin{align*}
&\bar{m}^h\circ \bar{y}^h=\mathscr{M}\quad\text{for every }h>0,\\
&\nabla_h(\bar{m}^h\circ \bar{y}^h)=(\nabla'\mathscr{M}|0)\quad\text{for every }h>0,\\
&u_{\bar{m}^h}\circ \bar{y}^h-\fint_{\Omega}u_{\bar{m}^h}\circ \bar{y}^h\,\md x\wk\mathscr{U}_{y,b,\mathscr{M}}\quad\text{weakly in }W^{1,2}(\omega;),\\
&\nabla_h (u_{\bar{m}^h}\circ \bar{y}^h)\wk (\nabla'\mathscr{U}_{y,b,\mathscr{M}}|\mathscr{V}_{y,b,\mathscr{M}})^T\quad\text{weakly in }L^2(\Omega;\R^3),
\end{align*}
and the following limsup inequality for the energy holds true:
$$\limsup_{h\to 0} E_h(\bar{y}^h,\bar{m}^h)\leq \mathcal{F}(y,b,\mathscr{M}).$$
 \end{theorem}
 \UUU A proof of the statement is in Section \ref{sec:proofs}
 below. \EEE

%%%%%%%%%%%%%%%%%%%%%%%%%%%%%%%%%%%%%%%%%%%%%%%%%%%%%%
\section{\UUU Proofs}\label{sec:proofs}

\UUU We collect in this section the proofs of the statements from
Section \ref{sec:gamma}. Within each subsection, notations are taken
from the corresponding statement. 

\subsection{Proof of Lemma \ref{lem:stray}}
 \EEE

We first observe that by \UUU the definition of the set of admissible
states $\mathcal{A}$ \EEE there holds 
\begin{equation}
\label{eq:ext-det-pos}
{\rm det}\,\widetilde{(\nabla'y|b)}\geq \ep\quad\text{on}\ \R^2.
\end{equation}
Additionally, for every $x'\in \R^2$ the matrix
$(\widetilde{(\nabla'y|b)}(x'))^{-1}(\widetilde{(\nabla'y|b)}(x'))^{-\UUU
  T \EEE}$ is symmetric. By \eqref{eq:ext-det-pos}, denoting by
$\lambda_i(x')$, $i=1,2,3$ the three eigenvalues of
$(\widetilde{(\nabla'y|b)}(x'))^{-1}(\widetilde{(\nabla'y|b)}(x'))^{-\UUU
  T \EEE}$ in increasing order, it follows that each of them is
different from zero for every $x'\in \R^2$. By the continuous
dependence of the eigenvalues of a matrix on the entries of the matrix
itself, and by the continuity of the map $x'\mapsto
(\widetilde{(\nabla'y|b)}(x'))^{-1}(\widetilde{(\nabla'y|b)}(x'))^{-\UUU
  T \EEE}$ (see \UUU again the definition of $\mathcal{A}$), \EEE we deduce that for every $i=1,2,3$ there exists a point $x^i\in \bar{\omega}$ such that
$$\min_{x\in \bar{\omega}}\lambda_i(x)=\lambda_i(x^i)>0.$$
Thus, recalling \eqref{eq:def-grad-ex}, we obtain
\begin{equation}
\label{eq:lambda}
\min_{i=1,2,3} \min_{x\in \R^2}\lambda_i(x)=\min_{i=1,2,3} \min\{1,\lambda_i(x^i)\}=:\lambda_{\rm eigen}>0.
\end{equation}
As a consequence of \eqref{eq:lambda}, the quadratic form 
$$Q(x,v):=(\widetilde{(\nabla'y|b)}(x'))^{-1}(\widetilde{(\nabla'y|b)}(x'))^{-\UUU T \EEE}v\cdot v\quad\text{for every}\,x'\in \R^2,\,v\in \R^3$$
satisfies
$$Q(x,v)\geq \lambda_{\rm eigen}|v|^2\quad\text{for every}\,x'\in \R^2,\,v\in \R^3.$$
The thesis is thus a direct consequence of the uniform ellipticity of
$Q$. \MMM \qed \EEE

\UUU
\subsection{Proof of Theorem \ref{prop:first-liminf}} \EEE
 
We subdivide the proof into three steps: in Step 1 we prove the compactness of sequences of deformations and magnetizations with equibounded energies. Step 2 is devoted to a characterization of the limiting stray field. Step 3 contains the proof of the liminf inequality.

\noindent\textbf{Step 1: Compactness}. In view of \eqref{eq:H1}, \eqref{eq:hp-Phi}, and \eqref{eq:uniform-en-estimate} we infer the existence of a constant $C$ such that
\begin{align}
&\label{eq:bd-lp-norm-grad}\|\nabla_h y^h\|_{W^{1,p}(\Omega;\mthree)}\leq C,\\
&\notag%\label{eq:bd-ls-neg-norm-det}
\Big\|\frac{1}{{\rm det}\,\nabla_h y^h}\Big\|_{L^q(\Omega)}\leq C,
\end{align}
for every $h>0$. By \eqref{eq:bd-lp-norm-grad}, and by the observation that
$$\|\nabla y^h\|_{L^p(\Omega;\mthree)}\leq \|\nh y^h\|_{L^p(\Omega;\mthree)},$$
 we deduce that there exists $y\in W^{2,p}(\Omega;\R^3)$ such that \eqref{eq:wk-null-av-def} is satisfied.
In particular, by \eqref{eq:bd-lp-norm-grad} we have $\partial_3 y=0$, thus $y$ can be identified with a map in $W^{2,p}(\omega;\R^3)$. As a further consequence of \eqref{eq:bd-lp-norm-grad}, we also find $b\in W^{1,p}(\omega;\R^3)$ and $d\in L^p(\Omega;\R^3)$ such that \eqref{eq:wk-conv-grad} and \eqref{eq:weak-d33} hold true. By \eqref{eq:wk-conv-grad}, the continuity of $\Phi$, and Fatou's lemma we obtain
\be{eq:lsc-inverse-det}\liminf_{h\to 0}\int_{\Omega}\Phi(\nabla_h y^h)\,\md x\geq \int_{\omega}\Phi(\nabla'y|b)\,\md x',\ee
which implies that ${\rm det}\,(\nabla'y|b)>0$ almost everywhere in
$\Omega$. Since $(\nabla'y|b)\in W^{1,p}(\omega;\mathbb{M}^{3\times
  3})\subset C^{0,\alpha}(\bar{\omega};\mathbb{M}^{3\times 3})$ for
$\alpha=(p-2)/p$, the argument in \cite[Theorem
3.1]{healey.kromer} yields \UUU $(\nabla'y|b)^{-1}\in
C^0(\bar{\omega};\mthree)$, $\text{det}(\nabla'y|b)\in
C^0(\bar{\Omega})$, and $\text{det}(\nabla'y|b)>\ep$ for some $\ep>0$. 

From convergences  \eqref{eq:wk-null-av-def}-\eqref{eq:wk-conv-grad} \EEE it follows in particular that 
\begin{equation}\label{eq:strong-det}{\rm det}\,\nabla_h y^h\to {\rm det}\,(\nabla' y|b)\quad\text{strongly in }C^0(\bar{\Omega}),\end{equation}
and hence 
\be{eq:lower-bd-det-h}{\rm det}\,\nabla_h y^h\geq \frac{\ep}{2}\quad\text{on }\bar{\Omega}\ee
for $h$ small. Properties \eqref{eq:uniform-en-estimate}  and
\eqref{eq:lower-bd-det-h} imply \UUU that \EEE
\begin{align}
&\nn\int_{\Omega}|(\nabla m^h)\circ y^h|^2\,\md x\leq \frac{2}{\ep}\int_{\Omega}|(\nabla m^h)\circ y^h|^2{\rm det}\,\nabla_h y^h\,\md x\\
&\label{eq:initial-bd-m}\quad=\frac{2}{h\ep}\int_{\Omega}|(\nabla m^h)\circ y^h|^2{\rm det}\,\nabla y^h\,\md x=\frac{2}{h\ep}\int_{\Omega^{y^h}}|\nabla m^h|^2\,d\xi\leq C.
\end{align}
In view of \UUU convergences \EEE \eqref{eq:wk-conv-grad} and \eqref{eq:strong-det}, there holds
\be{eq:unif-conv-inv}(\nabla_h y^h)^{-1}\to (\nabla'y|b)^{-1}\quad\text{strongly in }C^0(\bar{\Omega};\mathbb{M}^{3\times 3}),\ee
as well as
\be{eq:unif-conv}\nabla_h y^h\to (\nabla'y|b)\quad\text{strongly in }C^0(\bar{\Omega};\mathbb{M}^{3\times 3}).\ee
By combining \UUU bound \EEE \eqref{eq:initial-bd-m} with \UUU
convergence \EEE  \eqref{eq:unif-conv} we conclude that
\be{eq:bd-grad-mh-comp}
\int_{\Omega}|\nabla_h (m^h\circ y^h)|^2\,\md x\leq \int_{\Omega}|(\nabla m^h)\circ y^h|^2|\nabla_h y^h|^2\,\md x\leq C\int_{\Omega}|(\nabla m^h)\circ y^h|^2\,\md x\leq C.
\ee
In addition, by \eqref{eq:uniform-en-estimate} and by the saturation
constraint \UUU $|m|=1$   \EEE we deduce \UUU that \EEE
\be{eq:bd-mh-comp}
\int_{\Omega}|m^h\circ y^h|^2\,\md x\leq C.
\ee
Estimates \eqref{eq:bd-grad-mh-comp} and \eqref{eq:bd-mh-comp} yield
the existence of maps $\mathscr{M}\in W^{1,2}(\omega;\mathbb{S}^2)$
and $\eta\in L^2(\Omega;\R^3)$ such that \UUU convergences
\eqref{eq:m-wk} and \eqref{eq:m-wk-bis} hold, up to not relabelled
subsequences. \EEE 
In particular, there holds
$$(\nabla m^h)\circ y^h=(\nabla_h y^h)^{-T}\nabla_h (m^h\circ y^h)\wk (\nabla'y|b)^{-T}(\nabla'\mathscr{M}|\eta)\quad\text{weakly in }L^2(\Omega;\mthree),$$
and thus, by lower semicontinuity
\begin{align}
\label{eq:lsc-en-magnetization}
&\alpha\int_{\Omega}|(\nabla'y|b)^{-T}(\nabla'\mathscr{M}|\eta)|^2{\rm det}\,(\nabla'y|b)\,\md x
\leq \liminf_{h\to 0}\Big\{\frac{\alpha}{h}\int_{\Omega^y}|\nabla m|^2\Big\}.
\end{align}
The boundary conditions in \UUU the definition of $\mathcal{A}$ \EEE are a direct consequence of \eqref{eq:wk-conv-grad}. Thus, we conclude that $(y,b,\mathscr{M})\in~\mathcal{A}$.

Regarding the compactness of the stray field, we observe that by \eqref{eq:uniform-en-estimate}, \eqref{eq:lower-bd-det-h}, and \eqref{eq:unif-conv} there holds
\begin{align}
&\label{eq:bd-stray}\int_{\Omega}|\nabla_h (u_{m^h}\circ y^h)|^2\,\md
x\RRR\leq \UUU %\int_{\Omega}|(\nabla u_{m^h})\circ y^h|^2|\nabla_h y^h|^2\,\md x\leq C\int_{\Omega}|\nabla u_{m^h}\circ y^h|^2\,\md x\\
%&\nn\quad\leq C\int_{\Omega}|(\nabla u_{m^h})\circ y^h|^2{\rm det}\,\nabla_h y^h\,\md x
\frac{C}{h}\int_{\Omega^{y^h}}|\nabla u_{m^h}|^2\,d\xi\leq
\frac{C}{h}\int_{\R^3}|\nabla u_{m^h}|^2\,d\xi\leq C. \EEE
\end{align}
Therefore, by \UUU the \EEE Poincar\'e inequality we find $\mathscr{U}\in W^{1,2}(\omega;\R^3)$ and $\mathscr{V}\in L^2(\omega;\R^3)$ satisfying 
\begin{align*}
&%\label{eq:vh-one-almost} 
u_{m^h}\circ y^h-\fint_{\Omega}u_{m^h}\circ y^h\,\md x\wk\mathscr{U}\quad\text{weakly in }W^{1,2}(\omega),\\
&%\label{eq:vh-wk-almost} 
\nabla_h (u_{m^h}\circ y^h)\wk (\nabla'\mathscr{U}|\mathscr{V})^T\quad\text{weakly in }L^2(\Omega;\R^3).
\end{align*}

\noindent\textbf{Step 2: \UUU the Maxwell system}.
 In order to show that $\mathscr{U}=\mathscr{U}_{y,b,\mathscr{M}}$, $\int_{-\tfrac12}^{\tfrac12}\mathscr{V}dx_3=\mathscr{V}_{y,b,\mathscr{M}}$, and to pass to the limit in the magnetostatic energy, we observe that, since $u_{m^h}$ solves
 \be{eq:cond-umh}{\rm div}\,(-\mu_0 \nabla
 u_{m^h}+\chi_{\Omega^{y^h}}m^h)=0\quad \UUU \text{in} \ \R^3,\ee
 there holds
 \begin{align*}
 &\frac{\mu_0}{h}\int_{\UUU \R^3 \EEE}|\nabla u_{m^h}|^2\,d\xi=\frac{\mu_0}{h}\int_{\Omega^{y^h}}m^h\cdot \nabla u_{m^h}\,d\xi\\
 &\quad=\frac{\mu_0}{h}\int_{\Omega}(m^h\circ y^h)\cdot(\nabla u_{m^h})\circ y^h\, {\rm det}\nabla y^h\,\md x=\mu_0\int_{\Omega}(m^h\circ y^h)\cdot (\nabla_h y^h)^{-T}\nabla_h (u_{m^h}\circ y^h)\,{\rm det}\nabla_h y^h\,\md x.
 \end{align*}
 Therefore, by \eqref{eq:m-wk}, \eqref{eq:strong-det}, \eqref{eq:vh-wk}, and  \eqref{eq:unif-conv-inv} we conclude that
 \be{eq:conv-magn-en}
 \lim_{h\to 0}\frac{\mu_0}{h}\int_{\UUU \R^3\EEE}|\nabla u_{m^h}|^2\,d\xi=\mu_0\int_{\Omega}\mathscr{M}\cdot (\nabla'y|b)^{-T}(\nabla'\mathscr{U}|\mathscr{V})^T{\rm det}(\nabla'y|b)\,\md x.
 \ee

  We proceed now by passing to the limit into Maxwell's \UUU system. \EEE Denote by $\tilde{\Omega}$ the set 
  $$\tilde{\Omega}:=\mathbb{R}^2\times \big(-\tfrac12,\tfrac12\big),$$ 
  and consider the deformations
  \begin{equation}
  \label{eq:def-y-tilde-h}
  \tilde{y}^h(x):=\begin{cases} y^h(x)&\text{if}\,x\in\Omega\\
  (x',hx_3)&\text{if}\,x\in \tilde{\Omega}\setminus\Omega.
  \end{cases}
  \end{equation}
  In view of \eqref{eq:rescaled-bc} it follows that $\{\tilde{y}^h\}_h\subset W^{2,p}_{\rm loc}(\tilde{\Omega};\mathbb{R}^3)$. Let now $\varphi\in C^{\infty}_c(\tilde{\Omega})$. Choosing $\varphi\circ (\tilde{y}^h)^{-1}$ as a test function in \eqref{eq:cond-umh} we obtain that
 $$\frac1h \int_{\tilde{\Omega}^{\tilde{y}^h}}(\mu_0\nabla u_{m^h}-m^h)\cdot \nabla (\varphi\circ (\tilde{y}^h)^{-1})\,d\xi=0$$
 for every $h>0$. By performing a change of variables, the previous equation rewrites as
 \begin{equation}
 \label{eq:maxwell-h}
 \int_{\tilde{\Omega}}(\nh \tilde{y}^h)^{-1} [\mu_0(\nh \tilde{y}^h)^{-T}\nh (u_{m^h}\circ \tilde{y}^h)-\bar{m}^h\circ \tilde{y}^h]\cdot \nh \varphi\, {\rm det}\,(\nh \tilde{y}^h)\,\md x=0
 \end{equation}
 for every $h>0$ and $\varphi\in C^{\infty}_c(\tilde{\Omega})$, where 
 $$\bar{m}(\xi):=\begin{cases} m^h(\xi)&\text{if}\ \xi\in \Omega^{y^h}\\
 0&\text{otherwise in}\ \tilde{\Omega}^{\tilde{y}^h}.
 \end{cases}$$
 By \UUU the boundary conditions in $\mathcal{A}$,  convergences 
 \eqref{eq:unif-conv-inv} and  \eqref{eq:strong-det}, and by 
 definition \eqref{eq:def-y-tilde-h}, \EEE we deduce that
 \begin{align*}
 &%\label{eq:conv-inv-tilde}
 (\nh \tilde{y}^h)^{-1}\to \widetilde{(\nabla'y|b)}\quad\text{strongly in}\ C^0(\tilde{\Omega};\mthree),\\
 &%\label{eq:conv-det-tilde}
 {\rm det}(\nh \tilde{y}^h)^{-1}\to {\rm det}\widetilde{(\nabla'y|b)}\quad\text{strongly in}\ C^0(\tilde{\Omega}),
 \end{align*}
 where $\widetilde{(\nabla'y|b)}$ is the map defined in \eqref{eq:def-grad-ex}.
 Property \eqref{eq:m-wk} yields
 $$
 %\begin{equation}
% \label{eq:conv-ext-mh}
 \bar{m}^h\circ \tilde{y}^h\to \bar{\mathscr{M}}\quad\text{strongly in}\,L^2(\R^2), 
 %\end{equation}
 $$
 with $\bar{\mathscr{M}}$ as in \eqref{eq:def-m-bar}.
 Eventually, the same computations as in \eqref{eq:bd-stray} yield
 $$\int_{\tilde{\Omega}}|\nh (u_{m^h\circ \tilde{y}^h})|^2 dx\leq \frac{C}{h}\int_{\R^3}|\nabla u_{m^h}|^2 d\xi\leq C.$$
 Thus, by \eqref{eq:vh-one} and \eqref{eq:vh-wk} we deduce that there exist $\tilde{\mathscr{U}}\in W^{1,2}(\R^2)$ and $\tilde{\mathscr{V}}\in L^2(\tilde{\Omega})$ such that
 \begin{align*}
&u_{m^h}\circ \tilde{y}^h-\fint_{\Omega}u_{m^h}\circ \tilde{y}^h\,\md x\wk\tilde{\mathscr{U}}\quad\text{weakly in }W^{1,2}(\R^2),\\
&\nabla_h (u_{m^h}\circ \tilde{y}^h)\wk (\nabla'\tilde{\mathscr{U}}|\tilde{\mathscr{V}})^T\quad\text{weakly in }L^2(\tilde{\Omega};\R^3),
\end{align*}
with $\tilde{\mathscr{U}}=\mathscr{U}$ and  $\tilde{\mathscr{V}}=\mathscr{V}$ almost everywhere in $\Omega$.
 
Let now $\phi\in C^{\infty}_c(-\frac12,\frac12)$ and $\psi\in C^{\infty}_c(\R^2)$, and for every $h>0$ consider the function $\varphi^h(x):=\phi(hx_3)\psi(x')$ for every $x\in \R^2$. Choosing $\varphi^h$ as a test function in \eqref{eq:maxwell-h} for every $h>0$, and passing to the limit as $h\to 0$, we conclude that
 \begin{align*}
 &\int_{\R^2}\widetilde{(\nabla' y|b)}^{-1}[\mu_0 \widetilde{(\nabla' y|b)}^{-T}\Big(\nabla' \tilde{\mathscr{U}}|\int_{-\tfrac12}^{\tfrac12}\tilde{\mathscr{V}}\,\md x_3\Big)^T-\bar{\mathscr{M}}]\cdot (\nabla'\psi |0)^T\,{\rm det}\,\widetilde{(\nabla'y|b)}\phi(0)\,\md x\\
 &\quad + \int_{\R^2}\widetilde{(\nabla' y|b)}^{-1}[\mu_0 \widetilde{(\nabla' y|b)}^{-T}\Big(\nabla' \tilde{\mathscr{U}}|\int_{-\tfrac12}^{\tfrac12}\tilde{\mathscr{V}}\,\md x_3\Big)^T-\bar{\mathscr{M}}]\cdot (0 |\psi)^T\,{\rm det}\,\widetilde{(\nabla'y|b)}\phi'(0)\,\md x=0.%\label{eq:lim-max}
 \end{align*}
% for every $\phi\in C^{\infty}_c(-\frac12,\frac12)$ and $\psi\in
% C^{\infty}_c(\R^2)$. 
By the arbitrariness of \UUU $\phi\in C^{\infty}_c(-\frac12,\frac12)$
and $\psi\in C^{\infty}_c(\R^2)$ \EEE and by a density argument, we conclude that
%\begin{equation}
%\label{eq:Max1-wk} 
$$
\int_{\R^2}\widetilde{(\nabla' y|b)}^{-1}[\mu_0 \widetilde{(\nabla' y|b)}^{-T}\Big(\nabla' \tilde{\mathscr{U}}|\int_{-\tfrac12}^{\tfrac12}\tilde{\mathscr{V}}\,\md x_3\Big)^T-\bar{\mathscr{M}}]\cdot (\nabla'\psi |0)^T\,{\rm det}\,\widetilde{(\nabla'y|b)}\,\md x=0
$$%\end{equation}
for every $\psi\in W^{1,2}(\R^2)$, and
$$%\begin{equation}
%\label{eq:Max2-wk}
 \int_{\R^2}\widetilde{(\nabla' y|b)}^{-1}[\mu_0 \widetilde{(\nabla' y|b)}^{-T}\Big(\nabla' \tilde{\mathscr{U}}|\int_{-\tfrac12}^{\tfrac12}\tilde{\mathscr{V}}\,\md x_3\Big)^T-\bar{\mathscr{M}}]\cdot (0 |\psi)^T\,{\rm det}\,\widetilde{(\nabla'y|b)}\,\md x=0
$$%\end{equation}
for every $\psi\in L^{2}(\R^2)$. The identification $\mathscr{U}=\mathscr{U}_{y,b,\mathscr{M}}$ and $\int_{-\tfrac12}^{\tfrac12}\mathscr{V}\,\md x_3=\mathscr{V}_{y,b,\mathscr{M}}$ follows then by Lemma~\ref{lem:stray}.

\noindent\textbf{Step 3: Liminf inequality}. By \UUU convergences
\EEE\eqref{eq:wk-null-av-def}--\eqref{eq:weak-d33}, \UUU the liminf
inequalities 
\eqref{eq:lsc-inverse-det} and \eqref{eq:lsc-en-magnetization}, and
the  \UUU continuity  \EEE of $W$, we deduce that
\begin{align}
&\nn\liminf_{h\to 0} \Bigg\{\int_{\Omega}W(\nh y(x),m\circ y(x))dx+\frac{\alpha}{h}\int_{\Omega^y}|\nabla m(\xi)|^2d\xi+\int_{\Omega}|\nabla^2_h y(x)|^p\,\md x%\\
%\nonumber &\quad
+\int_{\Omega}\Phi(\nabla_h y(x))\,\md x\Bigg\}\\
&\nn\quad\geq \int_{\omega}W\big((\nabla' y|b),\,\mathscr{M}\big)\,\md x'+ \alpha\int_{\Omega}|(\nabla'y|b)^{-T}(\nabla'\mathscr{M}|\eta)|^2{\rm det}\,(\nabla'y|b)\,\md x\\
&\nn\qquad+\int_{\omega}\Bigg|\Bigg(\begin{array}{cc}(\nabla')^2y&\nabla'b\\(\nabla'b)^{\UUU T \EEE}&d\end{array}\Bigg)\Bigg|^p\,\md x'+\int_{\omega}\Phi(\nabla'y|b)\,\md x\\
&\nn\quad\geq\int_{\omega}W\big((\nabla' y|b),\,\mathscr{M}\big)\,\md x'+ \alpha\int_{\Omega}|(\nabla'y|b)^{-T}(\nabla'\mathscr{M}|0)|^2{\rm det}\,(\nabla'y|b)\,\md x\\
&\qquad+\int_{\omega}(|(\nabla')^2 y|^2+2|\nabla' b|^2)^{p/2}\,\md x'+\int_{\omega}\Phi(\nabla'y|b)\,\md x'. \label{eq:almost-liminf}
\end{align}
The liminf inequality \eqref{eq:liminf-inequality} follows by
combining \eqref{eq:conv-magn-en} with \eqref{eq:almost-liminf}, and
by \UUU recalling \EEE the characterization of the limiting stray field in Step 2. \MMM \qed \EEE
 
 \subsection{\UUU Proof of Theorem \ref{prop:limsup}} 
The statement follows by considering the following recovery sequences
$$y^h(x',x_3):=y(x')+hx_3 b(x')+f^h(x')-\fint_{\omega}f^h(x')\,\md x'$$
for almost every $x\in \Omega$, and
$$m^h(\xi):=\mathscr{M}\circ (y^h)^{-1}(\xi),$$
for almost every $\xi\in \Omega^{y^h}$, where $\mathscr{M}$ has been identified with a function defined on the infinite cylinder of basis $\omega$ and then has been extended to the whole $\R^3$. The convergence of the energies and the identification of the limiting stray field follow arguing as in the compactness argument. \MMM \qed \EEE
 
\subsection{\UUU Proof of Theorem \ref{thm:complete-gamma-conv}}
 The compactness and liminf inequality follow by \UUU Theorem
 \ref{prop:first-liminf} and by checking that property \eqref{eq:py}
 is preserved in the limit. \EEE The limsup inequality is obtained by observing that for $y$ satisfying \eqref{eq:py}, the maps $\bar{y}^h(x):=y(x')+hx_3 b(x')$ for every $x\in \R^3$ satisfy both \eqref{eq:average-inv} and \eqref{eq:ciarlet-necas}. The thesis follows by setting $$\bar{m}^h(\xi):=\mathscr{M}\circ (\bar{y}^h)^{-1}(\xi),$$
for almost every $\xi\in \Omega^{\bar{y}^h}$, where $\mathscr{M}$ has been identified with a function defined on the infinite cylinder of basis $\omega$ and then has been extended to the whole $\R^3$, and by arguing as in Proposition \ref{prop:limsup}. \MMM \qed \EEE

\section*{Acknowledgements}
\UUU We acknowledge support from the Austrian Science Fund (FWF)
projects F\,65, \RRR P\,29681, and \UUU  V\,662,  from the FWF-GA\v{C}R project
I\,4052/19-29646L, from the \RRR Vienna Science and Technology Fund (WWTF) projects \UUU MA14-009 and, \RRR partially supported also by Berndorf Privatstiftung and the City of Vienna, MA16-005, \UUU and from the OeAD-WTZ project CZ04/2019 (M\v{S}MT \v{C}R  8J19AT013).
% E.D. and U.S. acknowledge support from the Austrian Science Fund (FWF)
% through the special research program \emph{Taming complexity in
%   partial differential systems} (Grant SFB F65). The research of E.D.,
% M.K.,  and P.P. has been supported from the Austrian Science Fund
% (FWF) through the grant I 4052 N32 (GA\v{C}R project 19-29646L)
% \emph{Large Strain Challenges in Materials Science} and from BMBWF
% through the OeAD-WTZ project CZ04/2019 (M\v{S}MT \v{C}R  8J19AT013)
% \emph{Mathematical Frontiers in Large Strain Continuum Mechanics}. E.D
% has been also supported by the Austrian Science Fund (FWF) through the
% Elise-Richter grant V662 \emph{High contrast materials in plasticity
%   and magnetoelasticity}.  
\EEE
%
%P. Piovano acknowledges support from the Vienna Science and Technology Fund (WWTF), the City of Vienna, and Berndorf Privatstiftung under Project MA16-005, and the Austrian Science Fund (FWF) project P~29681.  

%\bibliographystyle{plain}
%\bibliography{ed}
\end{document}